\documentclass[reqno,11pt]{amsart}
\usepackage{ifthen, citesort}

\def\R{{\mathbb{R}}}
\def\C{{\mathbb{C}}}
\def\N{{\mathbb{N}}}
\def\ol#1{\overline{#1}}
\def\ul#1{\underline{#1}}
\def\fracp#1#2{\frac{\partial #1}{\partial #2}}
\def\theta{\vartheta}
\def\phi{\varphi}
\def\epsilon{\varepsilon}

\def\bbar{{\bar b}}
\def\cbar{{\bar c}}
\def\dbar{{\bar d}}
\def\ebar{{\bar e}}
\def\ibar{{\bar\imath}}
\def\jbar{{\bar\jmath}}
\def\kbar{{\bar k}}
\def\lbar{{\bar l}}
\def\qbar{{\bar q}}
\def\sbar{{\bar s}}
\def\dt{\frac{d}{dt}}
\def\g{g_{i\jbar}}
\def\tg{\tilde{g}_{i\jbar}}

\def\K{K\"ahler }
\def\A{Amp\`ere }
\def\H{H\"older }
\def\KR{K\"ahler-Ricci }
\def\KRF{K\"ahler-Ricci flow }
\def\nm{(n-1)}
\newcommand{\LOG}[1]{%
  \ifthenelse{#1 = 1}
  {\log\left(n^ns^{n-1}\right)}%
  {\log^{#1}\left(n^ns^{n-1}\right)}%
}
\DeclareMathOperator{\osc}{osc}

\newtheorem{theorem}{Theorem}[section]
\newtheorem{lem}[theorem]{Lemma}

\newtheorem{cor}[theorem]{Corollary}
\newtheorem{thm}[theorem]{Theorem}
\newtheorem{lemma}[theorem]{Lemma}
\newtheorem{corollary}[theorem]{Corollary}

\theoremstyle{definition}
\newtheorem{defn}[theorem]{Definition}
\newtheorem{definition}[theorem]{Definition}

\theoremstyle{remark}

\newtheorem{remark}[theorem]{Remark}
\newtheorem{notation}[theorem]{Notation}

\numberwithin{equation}{section}

\begin{document}

\title{Stability of Gradient K\"ahler-Ricci Solitons}

\author{Albert Chau}
\address{Harvard University, Department of Mathematics,
  One Oxford Street, Cambridge, MA 02138, USA}
\email{chau@math.harvard.edu}

\author{Oliver C.\ Schn\"urer}
\address{Max Planck Institute for Mathematics in the Sciences,
  Inselstr.\ 22-26, 04103 Leipzig, Germany}
\email{Oliver.Schnuerer@mis.mpg.de}

\renewcommand{\subjclassname}{%
  \textup{2000} Mathematics Subject Classification}

\date{April 2003, revised July 2003.}



\begin{abstract}
We study stability of non-compact gradient K\"ahler-Ricci flow solitons with
positive holomorphic bisectional curvature. Our main result is that any compactly
supported perturbation and appropriately decaying perturbations of the K\"ahler
potential of the soliton will converge to the original soliton under
K\"ahler-Ricci flow as time tends to infinity. To obtain this result, we
construct appropriate barriers and introduce an $L^p$-norm that decays for
these barriers with non-negative Ricci curvature.
\end{abstract}

\maketitle

\markboth{ALBERT CHAU AND OLIVER C. SCHN\"URER}
{STABILITY OF GRADIENT K\"AHLER-RICCI SOLITONS}

\section{Introduction}

In \cite{ham 3}, Hamilton introduced the Ricci flow to find Einstein metrics on a
compact Riemannian manifold by evolving an existing metric by its negative Ricci
curvature. By solving a complex Monge-Amp\`ere equation Yau \cite{yau} found
Einstein metrics on compact K\"ahler manifolds.  The results of Yau were later
re-established by Cao \cite{cao flow}, using the K\"ahler-Ricci flow
\begin{equation}\label{kr metric}
\dt g_{i\jbar}=-R_{i\jbar}
\end{equation}
to find Einstein metrics on a compact K\"ahler manifold. On $\C^n$
this equation corresponds to the evolution equation
\begin{equation}\label{kr potential}
\dt U=\log\det(U_{i\jbar})
\end{equation}
for the K\"ahler potential $U$ of the K\"ahler metric $g_{i\jbar}$. (For the
notations and conventions used in this paper, we refer to Section \ref{eq tr not
sec}.) In general Ricci and K\"ahler-Ricci flow may develop singularities well
before converging to an Einstein metric, and blow-up analysis gives rise to
complete non-compact solitons for these flow equations \cite{ham eternal,cao
resc}. Thus an understanding of the flow on non-compact manifolds is essential.
The general theory for the Ricci flow on non-compact manifolds was established in
a series of papers by Shi \cite{Shi} where in particular, parabolic maximum
principles on non-compact manifolds are established. For a survey article
concerning the Ricci flow and singularity analysis in particular see \cite{cao
chow}, a list of solitons for the K\"ahler-Ricci flow can be found in \cite{fik}.
Stability questions for the Ricci flow have been considered near compact Ricci
flat metrics and near complete metrics on $\R^2$ in \cite{gki,wu,hsu}.
Questions of stability and uniqueness in a limit as time tends to infinity of
(K\"ahler-)Ricci flow on compact manifolds have also been studied in many papers,
e.\ g.\ \cite{ham 3}. All these references represent only a small selection of
many articles concerning this subject. Another area of research is to rule out
special types of solitons after a blow-up by using geometric restrictions and to
prove uniqueness of those solitons.

In this paper, we focus on the non-compact complete gradient K\"ahler-Ricci
solitons found in \cite{Cao}. They are rotationally symmetric with positive
holomorphic bisectional curvature and their existence has been proved by solving
an ordinary differential equation. It turns out that these solutions are unique
(up to scaling and dilatations) in the class of rotationally symmetric gradient
solitons with positive holomorphic bisectional curvature. To learn more about
these solitons, it is desirable to know, whether they are stable under
appropriate perturbations. In this paper we answer this question in the
affirmative. We will show that the gradient solitons in \cite{Cao} on $\C^n$ are
stable under appropriately decaying perturbations of the K\"ahler potential.

 A further question in this direction is, whether there exist other solitons without
rotational symmetry. We wish to remark that both questions also seem to be
unsolved for the corresponding problems concerning strictly convex non-compact
hypersurfaces that arise as translating solitons for the mean curvature flow.

We wish to give a heuristic argument why we impose the condition that the
perturbation should decay at infinity. In our situation, we don't expect to get
stronger results than for the standard heat equation on $\R^n$. In this case,
however, we can take a bounded perturbation of the stationary solution $u=0$ that
satisfies for $t=0$
\begin{equation}\label{asy exa}
u((x_1,\, x_2,\,\ldots,\,x_n),\,t)\to
\begin{cases}
1 & \text{as~} x_1\to\infty,\\
-1 & \text{as~} x_1\to-\infty,
\end{cases}
\end{equation}
uniformly in $(x_2,\,\ldots,\,x_n)$. It follows directly
from the heat kernel representation of a solution
that \eqref{asy exa} remains true during the evolution,
i.\ e.\ for $t>0$. Similarly, we expect
that for general bounded perturbations of the potential in
our equation, the oscillation of a perturbation will
not tend to zero.  Of course in such special cases as the
one above one can show that the solution flattens out on
compact sets. This explains why it is natural to have decay
assumptions in our stability theorem. To simplify its
formulation, we give the following
\begin{definition}
A function $u_0:\C^n\to\R$ is called a C-potential, if it is rotationally
symmetric, $(u_0)_{i\jbar}$ is a K\"ahler metric with positive holomorphic
bisectional curvature, and gives rise to a gradient K\"ahler-Ricci soliton.
\end{definition}
We refer to Theorem 1 in \cite {Cao}, where C-potentials are shown to exist and
to be unique up to scaling and holomorphic transformations.  In particular, using
Notation \ref{fct conv}, it is shown that C-potentials are characterized by the
Equations \eqref{phi def 1} and \eqref{phi def 2}.

\begin{theorem}\label{stab thm}
Let $u_0$ be a C-potential in complex dimension $n\ge2$ and $\tilde u$ a smooth
perturbation such that
\begin{enumerate}
\item $\tilde{u}_{i\jbar}$ defines a complete K\"ahler metric on $\C^n$
equivalent to $(u_0)_{i\jbar}$ with bounded curvature. 
\item $u=\tilde u-u_0$
satisfies $\vert u(x)\vert\le K$ and
\begin{equation}\label{decay}
|u(x)|\le K\cdot \left(2\log\vert x\vert\right)^{-\alpha}
\quad\text{for~} |x|\ge 1
\end{equation}
for some $0<K$, $0<\alpha<1$.
\end{enumerate}

 Then with $u$ as initial
condition, \eqref{mpma2} has a long time smooth solution converging 
to $0$ as time tends to infinity.
\end{theorem}

In the special case of a compactly supported perturbation, we get

\begin{corollary}
Let $u_0$ be a C-potential in complex dimension $n\ge2$ and $\tilde u$ a smooth
perturbation such that
$\tilde{u}_{i\jbar}$ defines a complete K\"ahler metric on $\C^n$
equivalent to $(u_0)_{i\jbar}$ and $u=\tilde u-u_0$
is compactly supported.
 Then with $u$ as initial
condition, \eqref{mpma2} has a long time smooth solution converging 
to $0$ as time tends to infinity.
\end{corollary}

A geometric interpretation of the decay condition \eqref{decay}
is a follows. If the eigenvalues of the perturbed metric minus the
soliton metric $\tilde u_{i\jbar}-(u_0)_{i\jbar}$
with respect to the soliton metric 
$(u_0)_{i\jbar}$ decay like
$(2\log\vert x\vert)^{-2-\alpha}$ and $u$ tends
to $0$ at infinity, then \eqref{decay}
is fulfilled for an appropriate value of $K$. This is obtained by
integrating $u$ radially. Note that the barriers introduced
in \eqref{phi b} have the same decay in terms of the metric.

To give an overview over the method used here to prove Theorem \ref{stab thm}, we
describe our proof in words. It is convenient to transform our flow equation such
that the gradient solitons, we are interested in, become stationary solutions of
an equation, namely \eqref{mpma1}. In Section \ref{a priori}, we obtain smooth
longtime existence and get uniform estimates for the perturbation of a
C-potential. Then we want to apply the maximum principle to deduce that the
oscillation of the perturbation is strictly decreasing in time (or zero). Due to
the non-compactness of $\C^n$, however, we have to make sure that the supremum is
attained somewhere. We don't know how to prove this directly. Instead, we enclose
our perturbation from above and from below by radially symmetric barriers that
decay at infinity (as a function of $r=|x|$) and correspond to K\"ahler metrics
with positive holomorphic bisectional curvature. During the evolution, the upper
barrier stays positive, monotone in $r$, and rotationally symmetric. This ensures
that it attains its maximum at the origin and we can apply the strong maximum
principle to deduce that the oscillation, as a function of time $t$, is strictly
decreasing. However, this is not enough to show that the barrier converges to
zero. It might happen (and seems to be an interesting question, for which
equations it actually happens), that the fact that the perturbation tends to zero
at infinity is destroyed during the evolution as $t\to\infty$. This would imply
that the perturbation would converge to a positive constant as $t\to\infty$. For
the standard heat equation, however, it is quite easy to exclude this phenomenon
as the $L^2$-norm of a solution is non-increasing, so it remains finite during
the evolution, at least for $H^{1,2}$ initial data. The argument extends to any
smooth solution. In our situation, we can find a quantity (Lemma \ref{Lp lemma})
that is equivalent to the intrinsic $L^p$-norm of the perturbation and is also
decaying. This property relies heavily on the construction of special barriers
with positive holomorphic bisectional curvature (this is preserved during the
evolution \cite{Shi}). As the total (intrinsic) volume of our soliton is
infinite, this excludes the possibility that the perturbation tends to a positive
constant. Similar considerations apply to the lower barrier and thus the original
perturbation, enclosed in between these two barriers during the evolution,
converges to zero. Once $C^0$-convergence is established, smooth convergence
follows immediately from our a priori estimates.

In Section \ref{eq tr not sec}, we introduce the
notations that we will use throughout the paper and
transform the evolution equation \eqref{kr potential}
to other coordinate systems. We explain in Section
\ref{ste sec}, why we have shorttime existence of
solutions, and prove uniform a priori estimates
that guarantee longtime existence in Section \ref{a
priori}. We sketch how to construct barriers in
Section \ref{barrier construction} and refer to
the Appendices \ref{ODE sec} and \ref{phi inequ}
for details.
In the proof of Theorem \ref{stab thm}, we use the evident
Lemma \ref{monotonicity preserved} and give it's proof
in Appendix \ref{mon pres}.
Finally, Section \ref{to zero} contains the
proof, that our barrier converges to zero, the
crux of the proof of Theorem \ref{stab thm}.

The authors wish to thank J\"urgen Jost, Shing-Tung Yau,
the Alexander von Humboldt foundation
(Feodor Lynen Research Fellowship),
Harvard University, Cambridge, MA, U.S.A.,
and the Max Planck Institute for Mathematics
in the Sciences, Leipzig, Germany, for their support.
We thank Gerhard Huisken for his interest in the 
paper and the advice to reformulate the decay 
condition geometrically.

\section{Preliminaries and Transformations}
\label{eq tr not sec}

\subsection{Preliminaries}

\begin{notation}
We use indices to denote partial derivatives,
$$u_i=\fracp{}{z^i}u,\quad u_{i\jbar}=\fracp{^2}{z^iz^\jbar}u,
\quad\ldots.$$
If for a function $u:\C^n\to\R$, the matrix $(u_{i\jbar})$
is positive definite, we call $u$ a K\"ahler potential.
Then $(u_{i\jbar})$ is a K\"ahler metric and we denote
its inverse by $\left(u^{i\jbar}\right)$.
Lower case Latin indices range from $1$ to $n$. We use the Einstein
summation convention with a special convention for Latin capitals,
e.\ g.\
$$u^{i\jbar}w_{i\jbar}:=\sum\limits_{i=1}^n\sum\limits_{\jbar=1}^n
u^{i\jbar}w_{i\jbar},\quad z^Iu_I:=\sum\limits_{i=1}^nz^iu_i
+\sum\limits_{\jbar=1}^nz^\jbar u_\jbar.$$
We will use $\left(z^i\right)$ and
$\left(z^\jbar\right)$ to denote standard flat
coordinates on $\C^n$. Sometimes, it will be appropriate
to use standard Euclidean coordinates $\left(x^i\right)$.
The Laplace operator with respect to the metric $(u_{i\jbar})$
is defined by
$$\Delta_{u}w=u^{i\jbar}w_{i\jbar}.$$
In the estimates that follow, we will use $c$ to denote a
fixed positive constant that does not depend on time,
but may change its value from line to line.
Indices preceded by a comma, e.\ g.\ $u_{,i\jbar k}$,
indicate covariant differentiation with respect to
the background metric $(u_0)_{i\jbar}$ introduced in \eqref{mpma2}.
As usually, we use $R_{i\jbar}$ to denote the Ricci tensor,
$Rm$ for the Riemannian curvature tensor,
$\Vert\cdot\Vert$ to denote a (pointwise) norm with respect
to the induced metric, and $\nabla_g$ to indicate covariant
differentiation with respect to the metric $g_{i\jbar}$.
We do not use different notations for an the initial value
$u:\C^n\to\R$ and for the corresponding solution to
\KRF $u:\C^n\times[0,\infty)\to\R$.
\end{notation}

\begin{notation}
For a K\"ahler metric $(g_{i\jbar})$, we obtain Christoffel
symbols as follows
$$\Gamma^i_{jk}=g^{i\lbar}\fracp{g_{j\lbar}}{z^k},
\quad\Gamma^\ibar_{\jbar\kbar}=
\ol{\Gamma^i_{jk}},$$
and other components are identically zero.
Covariant differentiation is defined by
$$\omega^A_{,B}=\fracp{\omega^A}{z^B}+\Gamma^A_{BC}\omega^C,\quad
X_{A,B}=\fracp{X_A}{z^B}-\Gamma^C_{AB}X_C.$$
We can interchange covariant derivatives,
\begin{align*}
X_{c,ab}=&X_{c,ba},&
X_{\cbar,ab}=&X_{\cbar,ba},\\
X_{c,a\bbar}=&X_{c,\bbar a}-R_{a\bbar c\dbar}g^{\dbar e}X_e,&
X_{\cbar,a\bbar}=&X_{\cbar,\bbar a}
+R_{a\bbar d\cbar}g^{d\ebar}X_\ebar.
\end{align*}
In these formulae, the (holomorphic) Riemannian curvature
tensor appears, which is defined by
$$R_{i\jbar k\lbar}=
-\fracp{^2g_{i\jbar}}{z^k\partial z^\lbar}
+g^{p\qbar}\fracp{g_{i\qbar}}{z^k}\fracp{g_{p\jbar}}{z^\lbar}.$$
Contracting with respect to the metric yields the (holomorphic)
Ricci tensor
$$R_{i\jbar}=g^{k\lbar}R_{i\jbar k\lbar}
=-(\log\det (g_{k\lbar}))_{i\jbar}.$$
Finally, a K\"ahler manifold has positive biholomorphic
sectional curvature, if
$$R_{i\ibar j\jbar}>0.$$
Note, that we do not sum here.
\end{notation}

\begin{notation}\label{fct conv}
It will be convenient for the proofs to switch between different evolution
equations. If we assume in the following that a function $U$ fulfills \eqref{kr
potential}, we will also assume that the function $u$ is obtained from $U$ by
applying the transformations leading to \eqref{mpma2}, so it solves
\eqref{mpma2}. Similarly, we assume that $\tilde u$ fulfills \eqref{mpma1}. When
we consider a rotationally symmetric solution of \eqref{mpma1} depending on the
variable $s=\log|z|^2$, we denote this by $\hat u$. Analogous notations are used
for other functions solving \eqref{kr potential}.  Given a C-potential $u_0$, we
will denote by $U_0$ the corresponding solution to \eqref{kr potential}.

\end{notation}

\subsection{Transformations}

We will now fix a C-potential $u_0$.  Here and in the rest of the paper, all
evolution equations are defined on $\C^n$.

For further considerations, it will be convenient, to change coordinates such
that the evolution of a C-potential in time is as simple as possible. This can be
obtained by introducing
$$\tilde u(x,\,t):=U\left(e^{\frac12t} x,\,t\right)
+\tfrac12nt^2.$$ where $U$ is a solution to \eqref{kr potential}.  We
immediately get$\dt{\Tilde u}=\tfrac12 U_Iz^Ie^{\frac12t}+\dt U+nt$ and
$\det(\tilde u_{i\jbar})=e^{nt}\det (U_{i\jbar})$. Using \eqref{kr potential}, we
see that the evolution equation for $\tilde u$ is given by
\begin{equation}\label{mpma1}
\dt\tilde u=\log\det(\tilde u_{i\jbar})
+\tfrac12\left(z^i\tilde u_i
+z^\jbar \tilde u_\jbar\right),
\end{equation}
where we write again $z^i$ and $z^\jbar$ for
$e^{\tfrac12t}z^i$ and $e^{\tfrac12t}z^\jbar$,
respectively, i.\ e.\ we evaluate $\tilde u$ at
$\left(\left(z^i\right),\,
\left(z^\jbar\right),\,t\right)$.
The initial value is clearly
unchanged, $\tilde u(x,0)=U(x,0)$.

We now show that the C-potentials introduced in \cite{Cao} are in fact stationary
for \eqref{mpma1}. These potentials are characterized in \cite{Cao} by radially
symmetric functions $\hat{u}(s)$ in the variable $s=\log|z|^2$ for which the
following conditions are satisfied.  For $\hat u'(s)\equiv\phi(s)$ it is required
that $\phi(s)\to0$ for $s\to-\infty$ and $\phi$ fulfills (when normalized
appropriately) the ordinary differential equation
$$\phi^{n-1}\phi'e^\phi=e^{ns}.$$
By differentiating $\hat u\left(\log|z|^2\right) =\tilde u$ as
\begin{align*}
\tfrac12\left(z^i\tilde u_i+z^\jbar\tilde u_\jbar\right)
=&\hat u',\\
\tilde u_{i\jbar}=&\hat u''\frac{z_\ibar z_j}{\vert z\vert^4}
+\hat u'\frac1{\vert z\vert^2}\left(\delta_{i\jbar}
-\frac{z_\ibar z_j}{\vert z\vert^2}\right)
\end{align*}
 it is seen that a C-potential is a stationary solution to \eqref{mpma1}.   To show
 stability of \eqref{mpma1} at $u_0$, it will be convenient to write $\tilde u=u_0+u$. As $\tilde{u}$
and $u_0$ solve \eqref{mpma1} we get directly from the definition of $u$ its
evolution equation
\begin{equation}\label{mpma2}
\dt u=\log\frac{\det((u_0)_{i\jbar}+u_{i\jbar})}
{\det((u_0)_{i\jbar})}
+\tfrac12\left(z^iu_i+z^\jbar u_\jbar\right).
\end{equation}
The advantage of this evolution equation is that
$u_0$ is time-independent and it will turn out that
it allows to consider functions $u(x,\,t)$ that are
uniformly bounded in time.

\subsection{\H spaces}
 We now define the
 parabolic and elliptic \H spaces of a non-compact \K manifold $(M, \g)$.
 We will use these spaces to apply Schauder estimates in proving a priori
 estimates for \eqref{mpma2}.  These are parabolic versions of the elliptic
 spaces defined in \cite{CY1,TY,TY2}.

\begin{defn}\label{bddgeom} A \K manifold $(M,g_{i \bar\jmath})$
is said to have
bounded geometry of order $k+\alpha$, $k\in\N$, $0<\alpha<1$,
if there is a covering ${U_{i}}$ of M where,
\end{defn}

\begin{enumerate}
\item $U_{i}$ is holomorphically covered
by $\xi_{i} : V_{i} \rightarrow U_{i} $,
$V_i\subset\C^n$, for all i, where
$B_{r_{1}} (0) \subseteq V_{i} \subseteq
B_{r_{2}} (0)$ for $r_{1}>0$ and
$r_{2}<\infty$ independent of $i$.
\item Pulling the
metric back to $g^{*}_{a \bar{b}}$ on $V_{i}$ we get,
in flat coordinates on $V_{i}$,
 $(k_{1}\delta_{a \bar{b}} \leq g^{*}_{a \bar{b}}
\leq k_{2}\delta_{a \bar{b}})$, where $k_{1}>0$ and $k_{2}<\infty$ are
independent of $i$.  Also, for all $a$, $b$, we have
$\|g^{*}_{a \bar{b}}\|_{i,k +
\alpha} \leq C$ where the constant $C$ is independent of
$a$, $b$, and $i$ and
$\|\cdot\|_{i,k + \alpha}$ denotes the standard \H norm on $V_{i}$.
\end{enumerate}
A manifold is said to have bounded geometry of order infinity,
if it has bounded geometry of order $k+\alpha$ for every
$k\in\N$ and $0<\alpha<1$.

 Consider now the following norm defined for any smooth function $u$
on $M \times [0, T)$
\begin{equation}\label{holdernorm}
\indent \|u\|_{k + \alpha,
\frac{k}{2}+\frac{\alpha}{2}}:=\sup_{i}\{\|\xi_{i}^{*}u\|_{i,k + \alpha ,
\frac{k}{2}+\frac{\alpha}{2}}\}
\end{equation}
where $\xi_{i}^{*}u$ is the pull back of
$u$ to $V_{i}$ and $\|\cdot\|_{i,k +
\alpha, \frac{k}{2}+\frac{\alpha}{2}}$
is the standard parabolic \H norm on
$V_{i} \times [0, T)$ in the standard
coordinates on $V_{i}$.  We can now state
the following definition.

\begin{defn}
Let $(M, g_{i\bar\jmath})$ be a complete \K manifold of bounded geometry of order
$k+\alpha$. With respect to \eqref{holdernorm}, we define $C^{k + \alpha,
\frac{k}{2}+\frac{\alpha}{2}}(M \times [0, T))$ to be the closure of the set of
all smooth functions on $M \times [0, T)$ for which \eqref{holdernorm} is finite.
\end{defn}

When to time interval $[0, T)$ is
understood from the context, we will simply use
$C^{k + \alpha, \frac{k}{2}+\frac{\alpha}{2}}$ to denote the $C^{k + \alpha,
\frac{k}{2}+\frac{\alpha}{2}}(M \times [0, T))$.  Also, one can define the
elliptic \H spaces $C^{k + \alpha}$ in an obvious way.

\section{Short Time Existence}\label{ste sec}

We now establish the following general short time existence
result.
\begin{lem}\label{shorttime}
Let $(M, \g)$ be a complete non-compact \K manifold such that
$\Vert Rm\Vert\le c_0$ and $f:M\to\R$ is a smooth potential
of the Ricci tensor,
i.\ e.\ $R_{i\jbar}=-f_{i\jbar}$.
Then for some $T>0$ depending only on $c_0$, the following
initial value problem has a smooth solution $u(x, t)$ for $t\in
(0, T]$.
\begin{equation}\label{generalequation}
\begin{split}
 \frac{du}{dt}=&\log \frac{\det(\g+u_{i\jbar})}{\det(\g)}+f\\
 u(x, 0)&=0.
 \end{split}
\end{equation}
Moreover, for any $t\in(0, T]$, the \K metric $\g(x)+u_{i\jbar}(x, t)$ is
equivalent to $\g$ and has bounded geometry of order $\infty$ and $f(x)+(\log\det
(u_{i\jbar}))(x, t)$ is a potential for $R_{i\jbar}(x, t)$.
\end{lem}

This is the \K potential version of the following theorem of Shi \cite{Shi}.
\begin{thm}\label{shitheorem}
Let $(M,g_{i \bar\jmath})$ be a complete non-compact K\"{a}hler
manifold such that $\Vert Rm\Vert\leq c_0$.  Then for some constant
$T>0$ depending only on $c_0$,  there is a smooth short time solution
$\tilde{g}_{i \bar\jmath}(x,\,t)$  to the \K Ricci flow equation
\begin{equation}\label{krfequation}
\begin{split}
&\frac{d\tilde{g}_{i \bar\jmath}}{dt}= - \tilde{R}_{i \bar\jmath}
\\ &\tilde{g}_{i \bar\jmath}(x,0)=g_{i \bar\jmath}.
\end{split}
\end{equation}
for $t\in(0,\,T]$.  Moreover, for all $t\in(0, T ]$, $\tilde{g}_{i
\bar\jmath}(x,\,t)$ is a complete K\"{a}hler metric on $M$
equivalent to $g_{i \bar\jmath}$  and we have the following
estimates for the covariant derivatives of the curvature tensor of
$\tilde{g}_{i \bar\jmath}(x,\,t)$.
$$\left\Vert\nabla_{\tilde g} ^{m} \tilde{Rm}(x,\,t)\right\Vert^{2}
\leq C(n,m,c_0)(1/t)^{m}.$$
\end{thm}

\begin{proof}[Proof of Lemma \ref{shorttime}]
Under the hypothesis of the lemma, Theorem \ref{shitheorem}
guarantees a short time solution $\tg$ to the \KR flow
(\ref{krfequation}).  Using this solution, we solve the following
ordinary differential equation on $[0,\,T]$ for $x\in M$
\begin{align*}
 \frac{du}{dt}=&\log \frac{\det(\tg)}{\det(\g)}+f\\
 u(x, 0)=&0.
\end{align*}
for a smooth function $u(x, t)$.  It is then straight forward to
verify that we must have $\tg(x, t)=\g(x)+u_{i\jbar}$ and thus
$u(x, t)$ is a smooth solution to (\ref{generalequation}).  The
details of this verification can be found in \cite{ac}.    To
complete the proof of Lemma \ref{shorttime}, we need to show that
 for any $t\in(0,\,T]$ the \K metric
$\tg(x, t)=\g(x)+u_{i\jbar}(x, t)$ has bounded geometry of
order $\infty$. In \cite{TY2} the authors prove that on a non-%
compact K\"{a}hler manifold, one has bounded geometry of order $2+
\alpha$ provided one has bounded curvature, scalar curvature, and
gradient of scalar curvature.  Their proof can in fact be extended
to show that one has bounded geometry of infinite order provided
one has all covariant derivatives of curvature bounded.  Thus
since $\tg(x, t)$  has all covariant derivatives of its
curvature bounded, Theorem \ref{shitheorem}, we see that $\tg(x,
t)$ in fact has bounded geometry of order infinity.  This
completes the proof of Lemma \ref{shorttime}.
\end{proof}

\section{A Priori Estimates And Longtime Existence}
\label{a priori}

In this section we prove a priori estimates for solutions of \eqref{mpma2}.  We
follow the approach first used by Cao \cite{cao flow} which is to adapt the
elliptic estimates proved by Yau \cite{yau} and Aubin \cite{aubin} for the
elliptic complex Monge-Amp\`ere equation to the parabolic case. As we may
transform \eqref{mpma2} to an equation of the form \eqref{generalequation}, we
get short time existence. We may assume that we have a smooth solution $v\in
C^\infty(\C^n\times[0,\,T])$ to \eqref{mpma2}, that is, if $v$ is not smooth at
$t=0$, we use $t-\epsilon$, $1\gg\epsilon>0$, instead of $t$. Choosing $T$
smaller if necessary, we may also assume that $v(\cdot,\,t)$ gives rise to a
complete K\"ahler metric uniformly equivalent to $(u_0)_{i\jbar}$.

\subsection{Lower Order Estimates}

\begin{lemma}\label{zero order estimate}
A solution $v$ to \eqref{mpma2} satisfies
\begin{align*}
|v(\cdot, t)|_{C^{0}}\le&|v(\cdot, 0)|_{C^{0}}=:K_0\\
\intertext{and}
\left\vert\dt v(\cdot,\,t)\right\vert_{C^0}\le&
\left\vert\dt v(\cdot,\,0)\right\vert_{C^0}=:K_{\dt}
\end{align*}
for all $t\in[0,\,T]$.
\end{lemma}
\begin{proof}
This follows directly from the maximum principle in \cite{Ecker
Huisken}.
\end{proof}

\begin{lemma}
A solution $v$ to \eqref{mpma2} satisfies
$\left\vert z^Iv_I\right\vert\le c=:K_{z\nabla v}$
uniformly in $t$.
\end{lemma}
\begin{proof}
We estimate
\begin{align*}
\tfrac12\left\vert z^Iv_I\right\vert\le&
\left\vert\tfrac12 z^Iv_I-\dt v\right\vert
+\left\vert\dt v\right\vert
\le\left\vert\log\frac{\det(\tilde v_{i\jbar})}
{\det ((u_0)_{i\jbar})}\right\vert+K_{\dt}\\
=&\left\vert\log\frac{\det (V_{i\jbar})}
{\det ((U_0)_{i\jbar})}\right\vert+K_{\dt}
=\left\vert\dt V-\dt U_0\right\vert+K_{\dt}.
\end{align*}
As $(V-U_0)\left(e^{\frac12t}x,\,t\right)=(\tilde v-\tilde u_0)(x,\,t)
=v(x,\,t)$ is uniformly bounded in $C^0$, it suffices to prove
that
$$\left\vert\frac{d^2}{dt^2}V-\frac{d^2}{dt^2}U_0\right\vert\le c.$$
Then interpolation gives the claimed inequality. We differentiate
\eqref{kr potential} and obtain
\begin{align*}
\frac{d^2}{dt^2}V-\frac{d^2}{dt^2}U_0=&
\dt\log\det(V_{i\jbar})-\dt\log\det((U_0)_{i\jbar})\\
=&V^{i\jbar}\left(\dt V\right)_{i\jbar}-U_0^{i\jbar}
\left(\dt U_0\right)_{i\jbar}\\
=&V^{i\jbar}(\log\det(V_{k\lbar}))_{i\jbar}
-U_0^{i\jbar}(\log\det((U_0)_{k\lbar}))_{i\jbar}\\
=&-R_V+R_{U_0},
\end{align*}
where $R_V$ and $R_{U_0}$ are the scalar curvatures of the
metrics $V_{i\jbar}$ and $(U_0)_{i\jbar}$, respectively.
As $V$ and $U_0$ give rise to solutions of \KR flow, the
corresponding scalar curvatures are uniformly
bounded \cite{Shi}.
\end{proof}

We are not able to prove gradient estimates directly. Instead,
we have
\begin{lemma}\label{interpol lem}
Let $v$ be a solution to \eqref{mpma2}. Then there exists a
constant $K_{1+\alpha}$ that depends only on the C-soliton
such that
$$\Vert v(\cdot,\,t)\Vert_{1+\alpha}\le K_{1+\alpha}\cdot
\left(K_0+\left(n+\sup\limits_{\C^n}
\Delta_{u_0} v(\cdot,\,t)\right)\right)$$
for all $t\in[0,\,T]$.
\end{lemma}
\begin{proof}
We apply $L^p$-estimates \cite[Thm.\ 9.11]{gt} to
$0<n+\Delta_{u_0} v$ and obtain spatial $H^{2,p}$-bounds
for $v$. Then the Sobolev imbedding theorem implies the
result. Note that we only used $0<n+\Delta_{u_0}v$
and the $C^0$-bound.
\end{proof}

\subsection{Second Order Estimates}

 Consider the quantity
\begin{equation}
A=\log (n+\Delta v)-kv,
\end{equation}
where $\Delta v$ denotes the Laplacian of $v$ with respect
to $(u_0)_{i\jbar}(x)$ and the constant
$k\gg1$ is to be chosen later.
Clearly a bound on $|A|$ implies a bound on $|\Delta v|$.  We
will bound $A$ from above using the maximum principle.  The bound
from below will follow directly from some simple
inequalities.

\begin{lemma}\label{gradient A}
Assume that $\sup\limits_{\C^n}\Delta v\ge1$. Then
$$\frac{1}{n+\Delta v}\left(\Delta\left(z^I v_{I}\right)\right)\leq z^I
A_{I}+k\left(z^I v_{I}\right)+c$$
holds in $\Omega:=\left\{x\in\C^n:\Delta v(x)\ge
\tfrac12\sup\limits_{\C^n}\Delta v\right\}$.
\end{lemma}
\begin{proof}
We compute
\begin{equation}\label{1}
\begin{split}
\Delta\left(z^I v_{I}\right)
=&(u_0)^{l\bar{k}}\left(z^i v_{i}+z^{\jbar}
v_{\jbar}\right)_{l\bar{k}}\\
=&(u_0)^{l\bar{k}}\left(z^{i}_{,l}v_{i}
+z^i v_{,il}+z^{\jbar} v_{,\jbar l}\right)_{\bar{k}}\\
\le&(u_0)^{l\bar{k}}\left(z^{i}_{,l}v_{,i\bar{k}}
+z^{\jbar}_{,\bar{k}}v_{,\jbar l}+z^i
v_{,l\bar{k}i}+z^{\jbar} v_{,l\bar{k}\jbar }\right)
+c\cdot\Vert\nabla v\Vert_0
\end{split}
\end{equation}
and
\begin{equation}\label{2}
z^I A_{I}=\frac{1}{n+\Delta v}z^i\left((u_0)^{l\bar{k}}
v_{,l\bar{k}}\right)_{,i}
+\frac{1}{n+\Delta v}z^{\jbar}\left((u_0)^{l\bar{k}}
v_{,l\bar{k}}\right)_{,\jbar}-k z^Iv_I.
\end{equation}
Note that $z^i$, $z^\jbar$, and the Riemannian curvature tensor
induced by $u_0$ are bounded with respect to the metric
$((u_0)_{i\jbar})$, see Remark \ref{bounded rem}. This
allows to estimate the terms obtained by interchanging
the order of covariant differentiation.

We combine (\ref{1}) and (\ref{2}), and get in $\Omega$
\begin{equation}\label{3}
\begin{split}
\frac{1}{n+\Delta v}&\Delta\left(z^I v_{I}\right)
-z^I A_{I}\\
\leq &k z^I v_{I} + c \frac{\Vert\nabla v\Vert_0}{n+\Delta v}
+\frac{1}{n+\Delta v}
(u_0)^{l\bar{k}}\left(z^{i}_{,l}v_{i\bar{k}}+z^{\jbar}_{,\bar{k}}
v_{\jbar l}\right)\\
\leq &k z^I v_{I}+c\frac{\Delta v}{n+\Delta v} +
c\frac{1}{n+\Delta v}.
\end{split}
\end{equation}
Here we have used that at a fixed point,
we can always choose holomorphic coordinates
such that $(u_0)^{l\bar{k}}=\delta^{l \kbar}$ and $v_{i\bar{k}}=0$ for
$i\neq k$, so we get
$$(u_0)^{l\bar{k}}\left(z^{i}_{,l}v_{i\bar{k}}+z^{\jbar}_{,\bar{k}}
v_{\jbar
l}\right)\leq c\Delta v+ c$$
and deduce the second inequality in \eqref{3}.
 In such coordinates, the terms $1+v_{i\ibar}$
are positive for each $i$ and are simply the eigenvalues of the
tensor $(u_0)_{i\jbar}+v_{i\jbar}$ with respect to metric
$(u_0)_{i\jbar}=\delta_{i\jbar}$.
Finally, in passing to the last line of
(\ref{3}), we have used Lemma \ref{interpol lem}.
\end{proof}

We are now in a position to prove the following
\begin{lemma}\label{laplace estimate}
There is a constant $K_{2}>0$  such that
$\vert\Delta v(x,\,t)\vert\le K_{2}$ for all
$(x,\,t)\in\C^n\times\in[0, T]$.
\end{lemma}
\begin{proof}
For the proof of the upper bound for $\Delta v$, we will only consider
those $(x,\,t)\in
\C^n\times[0,\,T]$ such that $\sup\limits_{\C^n\times\{t\}}
\Delta v>1$ and
$\Delta v(x,\,t)>\tfrac12\sup\limits_{\C^n}\Delta v(\cdot,\,t)$.
Thus, we can use Lemmata \ref{interpol lem} and \ref{gradient A}.
We compute the evolution equation for $A$, interchange fourth
covariant derivatives and use \cite[p.\ 264]{aubin} to estimate
third derivatives
\begin{align*}
\frac{dA}{dt}-\tilde\Delta A\le&
\frac1{n+\Delta v}\frac{d\Delta v}{dt}-k\frac{dv}{dt}
-(k-c)\tilde v^{i\jbar}(u_0)_{i\jbar}\\
&-\frac1{n+\Delta v}
\left(\Delta\frac{dv}{dt}-\tfrac12\Delta\left(z^Iv_I\right)\right)+nk,
\end{align*}
where $\tilde\Delta$ denotes the Laplacian in the metric
$\tilde v_{i\jbar}(x, t)=(u_0)_{i\jbar}(x)+v_{i\jbar}(x, t)$ with
inverse $\tilde v^{i \jbar}$.
We use the geometric-arithmetic means inequality and
$$\left(\sum\limits_{i=1}^n\frac1{\lambda_i}\right)^{n-1}
\ge\tfrac1n\sum\limits_{i=1}^n\lambda_i\cdot\prod\limits_{i=1}^n
\frac1{\lambda_i},\quad\lambda_i>0,$$
which is proved easily as we may assume
that $1=\lambda_1\le\lambda_2\le\ldots
\le\lambda_n$, to obtain
\begin{align*}
(n+\Delta v)\geq& \left[\frac{\det(\tilde v_{i \jbar})}{\det
((u_0)_{i\jbar})}\right]^{\frac{1}{n}} =
e^{\frac{1}{n}\left(-\frac12z^I v_{I}+\frac{dv}{dt}\right)},\\
\tilde v^{i\jbar}(u_0)_{i\jbar}\ge&\frac1{c(n)}
\left[(n+\Delta v)\cdot e^{\frac12 z^Iv_I-\frac{dv}{dt}}
\right]^{\frac1{n-1}}.
\end{align*}
We apply Lemma \ref{gradient A} and estimate
\begin{align*}
\frac{dA}{dt}-\tilde\Delta A \leq &\frac{1}{n+\Delta v}
\frac{d \Delta v}{dt} -k \frac{dv}{dt}-(k-c) \tilde v^{i\jbar}
(u_0)_{i\jbar}\\
 &- \frac{1}{n+\Delta v} \Delta \frac{dv}{dt}
+\tfrac12 z^I A_{I}+\tfrac12 k z^I v_{I} + nk+c\\
 \leq &-(k-c) \tilde v^{i\jbar}(u_0)_{i\jbar}
+\tfrac12 z^I A_{I}
+\tfrac12 k z^I v_{I} +nk-k \frac{dv}{dt}+c\\
 \leq & -\frac{k-c}{c(n)}e^{\frac{1}{n-1}\left(\frac12z^I
v_{I}-\frac{dv}{dt}\right)}
(n+\Delta v)^{\frac{1}{n-1}}\\
&+\tfrac12 z^I A_{I}+\tfrac12 k z^I v_{I}+nk -k\frac{dv}{dt}+c.
\end{align*}
Fixing $k\gg1$ so large that $k-c$ is bounded below by some
positive constant, the maximum principle can be applied to
the evolution equation
$$\dt A-\tilde\Delta A\le-\tfrac1ce^{\frac{A-c}{n-1}}
+\tfrac12 z^IA_I+c,$$
implying the upper bound.

To prove a lower bound for $\Delta v$, we use coordinates as
in Lemma \ref{gradient A}. Our lower order estimates imply
that $\prod\limits_{i=1}^n(1+v_{i\ibar})$ is bounded below
by a positive constant. The function $v$ gives rise to a
\K metric. So all the factors are positive.
As we have seen that $v_{i\ibar}$
is uniformly bounded above for each $i$, the lower bound
follows.
\end{proof}

\begin{cor}
There is a constant $K_{1}>0$ depending only on $K_{0}$ and
$K_{2}$ such that $\Vert v(\cdot,t)\Vert_1<K_{1}$
for all $t\in[0, T]$.
\end{cor}
\begin{proof}
Use Lemmata \ref{zero order estimate}, \ref{interpol lem}
and \ref{laplace estimate}.
\end{proof}

\begin{cor}\label{equivalence}
The metric $w_{i\jbar}(x, t)=v_{i\jbar}(x,
t)+(u_0)_{i\jbar}(x)$ is equivalent to $(u_0)_{i\jbar}(x)$ for all
$t\in[0, T]$. Moreover, the equivalence factor depends only on
$K_{0}$, $K_{1}$ and $K_{2}$.
\end{cor}
\begin{proof}
This follows from the proof of Lemma \ref{laplace estimate}.
\end{proof}

\subsection{Higher Order Estimates and Long Time Existence}

Consider the quantities
\begin{align*}
Q_3= \tilde v^{i\jbar}\tilde v^{k\lbar}\tilde v^{r\sbar}
  v_{,i\lbar r}v_{,\jbar k\sbar},\\
\intertext{and}\\
Q_4 = \tilde v^{i\jbar}\tilde v^{k\lbar}\tilde v^{r\sbar}
  \tilde v^{a\bbar} v_{,i\lbar r\bbar}v_{,\jbar k\sbar a},\\
\end{align*}
where the covariant differentiation is with respect to
$(u_0)_{i\jbar}$ and $\tilde v^{i \jbar}$ represents the
inverse of the time dependent metric $\tilde v_{i\jbar}(x, t)
=(u_0)_{i\jbar}(x)+v_{i\jbar}(x, t)$.  By the
previous section, this norm is equivalent to that using
$(u_0)_{i\jbar}$.

\begin{lem}\label{higher}
There are constants $K_{3}, K_{4}>0$
depending only on $K_{0}$, $K_{1}$, $K_{2}$
such that $|Q_3|_{C^0}<K_{3}$ and $|Q_4|_{C^0}<K_{4}$
for all $t\in[0, T]$.
\end{lem}
\begin{proof}
The above estimates are known in the special
case that $v$ is a solution to
(\ref{generalequation}) and have appeared
in several places in various equivalent
forms. We describe some of these briefly.
Calabi first estimated $|Q_3|_{C^0}$ for
the elliptic Monge-\A equation on a
compact manifold. This estimate was later
used by Aubin \cite{aubin} and by Yau
\cite{yau} in proving the Calabi conjecture.
Calabi's estimate was applied directly
by Cao \cite{cao flow} to (\ref{generalequation})
and later by Shi \cite{Shi} to (\ref{kr metric}).
In \cite{Shi} Shi goes further to estimate an
appropriate second derivative of the
solution to (\ref{kr metric}) and observes
that this is equivalent to estimating
the curvature tensor of the evolving metric.
An equivalent estimate can be found
in \cite{ac} where $|Q_4|_{C^0}$ is
estimated for (\ref{generalequation}).

In our case that $v$ is a solution to (\ref{mpma2}), we point out
that it is straight forward to adapt the arguments of the authors
cited above to our case.
\end{proof}
Notice that while the estimates in Lemma \ref{higher} follow rather painlessly
from the corresponding estimates for (\ref{generalequation}), such is not the
case for our laplacian estimate in Lemma \ref{laplace estimate}. The difference
is that in Lemma \ref{higher} we have already estimated all derivatives of lower
order and second derivatives of the form $v_{i\jbar}$, 
while in the case of Lemma \ref{laplace estimate} do not have gradient
estimates a priori. This does not cause a problem in (\ref{generalequation})
while in our case it does.

Note that the a priori estimates obtained so far imply that
for any $t\in[0,T]$, $n+\Delta_{u_0}v(\cdot,t)\in C^\alpha$
with uniform bounds. Thus elliptic Schauder theory implies
that $v(\cdot,t)\in C^{2+\alpha}$. Differentiating \eqref{mpma2}
yields $v_{i\jbar}\in C^{\alpha,\frac{\alpha}2}$ and $v\in C^{2+\alpha,
1+\frac\alpha2}$ with uniform bounds.

\begin{lem}
Let $v$ be a solution to \eqref{mpma2}
and let $C^{k+ \alpha,
\frac k2+\frac{\alpha}2}$ be the \H spaces on $\C^n$ relative
to the metric $\tilde v_{i\jbar}$.  Then for
every $k$, $v$ is bounded in $C^{k + \alpha,
\frac k2+\frac{\alpha}2}$ independent of $t$.
\end{lem}
\begin{proof}
We prove the respective result for H\"older spaces with respect
to the background metric. The corresponding results in these 
H\"older spaces imply the claimed estimates.
 Consider an arbitrary coordinate
neighborhood $V_\beta$ with coordinates $\left(z^i\right)$
as in Definition \ref{bddgeom}.
Differentiating (\ref{mpma2}) with respect to
$z^i$ in this coordinates and rearranging terms gives
\begin{equation}\label{bootstrapequation}
\begin{split}
\dt v_i=&\tilde v^{r\sbar}v_{,r\sbar i}
-\tfrac12\left(\hat z^I u_I\right)_i\\
=&\tilde v^{r\sbar}v_{,ir\sbar}
-\tfrac12\hat z^Iv_{,iI}
+\left(\tilde v^{r\sbar}
R_{i\sbar r\dbar}\tilde u_0^{e\dbar} v_{e}-\tfrac12 \hat z^I_{,i}
u_I\right),
\end{split}
\end{equation}
 where $\hat z^i$ are the local components of
the global vector field $z^i$.
For covariant differentiation and the curvature
tensor, we use the background metric $(u_0)_{i\jbar}$.
 We view (\ref{bootstrapequation}) as a
parabolic equation for $v_i(x, t)$ on the
coordinate domain $V_\beta \times [0, \infty)$
with the third term on the right-hand
side considered as a single inhomogeneous
term. In what follows, all bounds
stated will be independent of $\beta$ and $t$.
It is readily seen that our estimates
from above provide us with
a $C^{\alpha, \frac{\alpha}{2}}$
bound for the coefficients and terms of
(\ref{bootstrapequation}). We may then
apply standard parabolic Schauder
estimates to obtain a $C^{2 + \alpha,
1+\frac{\alpha}{2}}$ bound for $v_i(x, t)$
in an interior domain of $V_\beta$.  A
standard bootstrapping argument \cite{cao flow}
combined with the fact that the metric
$\tilde v_{i\jbar}$ has bounded geometry of order
$\infty$ then allows us to obtain a
$C^{k + \alpha, \frac{k}{2} + \frac{\alpha}{2}}$
bound on $v(x, t)$ in $V_\beta$ for
all $k$. The lemma now follows readily from Definition \ref{bddgeom}.
\end{proof}

\begin{cor}
The solution $v$ is smooth and exists for all time.
Moreover the metric $\tilde v_{i\jbar}(x, t)=(u_0)_{i\jbar}(x)
+v_{i\jbar}(x,\,t)$ remains equivalent to $(u_0)_{i\jbar}(x)$
uniformly over all $t$ and the curvature of
$\tilde v_{i\jbar}(x, t)$
remains bounded on $\C^{n}$ independent of $t$.
\end{cor}
\begin{proof}
By our a priori estimates, it is straight forward to see that the
curvature of the metric $\tilde v_{i\jbar}(x, t)$ stays uniformly
bounded on $[0, T]$.
Corollary \ref{equivalence} implies that the metrics stay uniformly
equivalent.
Thus to prove the corollary it suffices to
prove the assertion of long time existence.  Moreover, long time
existence for $v$ follows from long time existence for
\eqref{kr potential}
with initial condition $\tilde v(x, 0)=u_{0}(x)+v(x, 0)$. Begin by
assuming that $T$ is the maximal time up to which we have a smooth
solution. Choosing a time $T'<T$ arbitrarily close to $T$ and
applying Lemma \ref{shorttime} to the metric $\tilde v_{i\jbar}(x,
T')$, we may extend $\tilde v$ past $T$ as a solution to
\eqref{kr potential} thus
arriving at a contradiction and thus proving the corollary.
\end{proof}

\section{Barrier Construction}
\label{barrier construction}

Before we can construct a barrier, we have to determine
the precise asymptotic behavior of our soliton. According
to \cite{Cao}, we may assume that the function
$\phi(s)=\hat u'(s)$ fulfilling
\begin{equation}\label{phi def 1}
\phi^{n-1}\phi'e^\phi=e^{ns},
\end{equation}
and
\begin{equation}\label{phi def 2}
\phi(s)\to0\quad\text{for~} s\to-\infty,
\end{equation}
where $s=\log|z|^2$, gives rise to our soliton.
The second condition is required to obtain a
smooth solution at the origin.
We derive in Appendix \ref{ODE sec}
the following expansions for $\phi$ and its derivatives
at infinity
\begin{align*}
\phi=&ns+o(s),\displaybreak[1]\\
\phi'=&n+o(1),\displaybreak[1]\\
\phi''=&\frac{n-1}{s^2}+o\left(\frac1{s^2}\right),\\
\intertext{and}
\phi'''=&-2\frac{n-1}{s^3}+o\left(\frac1{s^3}\right).
\end{align*}

In the following, we construct barriers in the case
$n\ge2$.
Now, we assume that our perturbation $u(x,0)$ of the initial
value is such that
$$|u(x,0)|\le K\cdot\min\left\{1,s^{-\alpha}\right\},
\quad\text{where~}s=2\log|x|,\quad 0<\alpha<1.$$
For our barrier we make the ansatz
\begin{equation}\label{phi b}
\phi_b(s)=\phi(s)\mp K s^{-1-\alpha}\alpha(2R)^\alpha
\psi\left(\frac sR\right)
\end{equation}
with $\phi$ as above, that corresponds to the barrier
$$\hat b(s)=\hat u_0(s)\pm\int\limits_s^\infty K\sigma^{-1-\alpha}
\alpha(2R)^\alpha\psi\left(\frac \sigma R\right)d\sigma.$$
Here $\psi$ is a smooth monotone function such that
$$\psi(s)=\begin{cases}
0&\text{if~}s\le1,\\
1&\text{if~}s\ge2.\\
\end{cases}$$
Assume from now on that $R\ge\tfrac12$.
It is straight-forward to check that $\hat b(s)$ lies
above/below our perturbed initial value.

To prove that for $R\gg1$ fixed sufficiently large
\begin{equation}\label{5 cond}
\phi_b>0,\quad
\phi_b'>0,\quad
\phi_b-\phi'_b>0,\quad
\phi_b'^2-\phi_b\phi_b''>0,\quad
\phi_b''^2-\phi_b'\phi_b'''>0
\end{equation}
is again a technical calculation, we refer to
Appendix \ref{phi inequ}.

Note that it is essentially the integrability condition
for $\phi_b$ and not \eqref{5 cond} that determines
the possible exponents in the decay
condition.

For $n=1$, our method does not seem to work.
In this case, $\phi(s)$ is even explicitly known
to be $\log\left(1+e^s\right)$, but
$$\phi''^2-\phi'\phi'''=e^{-s}+O\left(e^{-2s}\right)$$
seems to exclude such a barrier construction.
For results concerning longtime behavior of
solutions to Ricci flow in the corresponding real
dimension $2$, we refer to \cite{hsu,wu}.

\section{Convergence to Zero}\label{to zero}

In this section we prove Theorem \ref{stab thm}. We use the radially symmetric
decaying barriers constructed in Section \ref{barrier construction} to enclose
our initial perturbation from above and from below.  By smoothly evolving our
barriers and perturbed initial value using (\ref{mpma2}) for all time, the
maximum principle of \cite{Ecker Huisken} implies that our perturbation will
converge to zero provided such is true of our barriers.  In particular, the
perturbed soliton converges back to the original soliton as $t\to\infty$.  We
will only show that the upper barrier converges back to the original soliton.
Studying the behavior during \KRF is simpler for the barriers as they are
rotationally symmetric and decaying in $|z|$.

\begin{lem}\label{convergencetozero}
Let $b$ be the upper (lower) barrier constructed in Section
\ref{barrier construction}. Then (\ref{mpma2}) with
initial condition $b$ has a long time smooth solution, which
we also denote by $b$, which converges to zero as
$t\to\infty$ in the $C^0$ norm.
\end{lem}

 The proof is divided into several steps.  We sketch the proof
for the case of the upper barrier and note that the
 case of the lower barrier is similar.
Part of the argument is a modification of the convergence
proof in \cite{os}.
We first show that the
condition that $b$ initially decays monotonely in $|z|$ is
 preserved for all time, so we get
especially $b(0,\,t)\ge b(x,\,t)$ for all
 $(x,\,t)\in\C^n\times[0,\,\infty)$. We do this in Lemma
\ref{monotonicity preserved}.  The strong maximum principle then
 guarantees that $\sup_{\C^n}b$ is strictly decreasing
 in $t$.  In fact, we claim that $b$ must converge to a constant.  This
can be seen as follows. In view of our a priori estimates, we can find
for every sequence
$t_n\to\infty$ a subsequence, again denoted by $t_n$,
such that the maps
$$\C^n\times[-t_n,\infty)\ni(x,\,t)\mapsto b(x,\,t+t_n)$$
converge locally uniformly in any $C^k$-norm to a smooth function
$b^\infty(x,\,t)$ satisfying the evolution equation (\ref{mpma2}) everywhere in
$\C^n\times\R$. Moreover, since the oscillation of $b$ decreases strictly in time
by the strong maximum principle \cite{aw}, it must converge to some nonnegative
constant. In other words, the limit solution $b^\infty(x,\,t)$ has nonnegative
oscillation which is constant in time. But it is easy to see that the rotational
symmetry and decay condition on $b(x, t)$ also holds for $b^{\infty}(x, t)$ and
thus by the strong maximum principle, the oscillation of $b^{\infty}(x, t)$
cannot be a positive constant. Thus $b^{\infty}(x, t)$ is constant in space. The
monotonicity of $b(0,\,t)$ shows that this constant is independent of the chosen
subsequence and hence $b$ actually converges to a constant. In Corollary \ref{Lp
corollary} we show that during the evolution the $L^p$-norm, for some $p\ge2$, of
$b$ is dominated by its value at $t=0$. We compute the $L^p$-norm with respect to
an evolving volume form which stays uniformly equivalent to the volume form for
the initial soliton metric, thus the integral of any positive constant over
$\C^n$ with respect to this volume form for fixed $t$ is infinite. By Remark
\ref{ip fin}, the $L^p$-norm is finite for $t=0$ provided $p>\tfrac{n+1}\alpha$
with $\alpha$ as in Section \ref{barrier construction}. So $b$ has to converge
uniformly to zero on compact subsets of $\C^n$ as $t\to\infty$. Note that the
monotonicity in $|z|$ is preserved during the evolution and when we extract
subsequences. Moreover, $b(0,t)$ is decreasing in $t$. So $b(0,\,t)$ has to
converge to zero and it follows that the perturbed soliton converges back to the
original soliton in $\C^n$.

\begin{lemma}\label{monotonicity preserved}
Let $b$ be the upper barrier constructed in Section \ref{barrier construction}.
Then $b$ stays rotationally symmetric and the property that $b$ decays in $|z|$
is preserved during the evolution of $b$ by \eqref{mpma2}.
\end{lemma}
It is quite evident that this lemma is true. Thus,
we defer it's proof to Appendix \ref{mon pres}.

\begin{lemma}\label{Lp lemma}
Let $u_0$ be a C-potential and $b$ a barrier as constructed in Section
\ref{barrier construction}. Then there exists a metric $a_{i\jbar}$, uniformly
equivalent to $(u_0)_{i\jbar}$ and $(u_0)_{i\jbar}+b_{i\jbar}$, such that for $b$
evolving according to (\ref{mpma2}), we have
$$\int\limits_{\C^n}\vert b(t)\vert^p\det(a_{i\jbar}(t))
\equiv I_p(t)\le I_p(0)$$
for $p\ge2$ and $t\ge0$.
\end{lemma}
\begin{proof}
Interpolating between the two determinants in (\ref{mpma2})
and using upper indices to denote inverses, we get
\begin{align*}
\dt b=&\log\det((u_0)_{i\jbar}+b_{i\jbar})-\log\det((u_0)_{i\jbar})
+\tfrac12 z^Ib_I\\
=&\int\limits_0^1((u_0)_{\cdot\cdot}+\tau b_{\cdot\cdot})^{i\jbar}
d\tau b_{i\jbar}+\tfrac12 z^Ib_I\\
\equiv&a^{i\jbar}b_{i\jbar}+\tfrac12 z^Ib_I.
\end{align*}
Now we define $(a_{i\jbar})$ to be the inverse of $\left(a^{i\jbar}\right)$. By
definition, $(a_{i\jbar})$ is uniformly equivalent to $(u_0)_{i\jbar}$ and
$(u_0+b)_{i\jbar}$ as these two metrics stay uniformly equivalent during the
evolution. For showing the definiteness of terms like $\dt (u_0)_{i\jbar}$ and
$\dt ((u_0)_{i\jbar}+b_{i\jbar})$, it will be convenient to substitute so that we
almost come back to the original evolution equation (\ref{kr potential}). Set
$\ul b(x,\,t):=b\left(e^{-\frac12t}x,\,t\right)$, $\ul
u_0(x,t):=u_0\left(e^{-\frac12t}x,t\right)$. This implies that $(\ul u_0+\ul
b)(x,\,t)= B(x,\,t)+\tfrac12nt^2$, where $B=B(u_0+b)$ is as in \eqref{kr
potential}. As the metric $B_{i\jbar}$ has positive holomorphic bisectional
curvature for $t=0$ (Appendix C), this is preserved during the evolution
\cite{Shi}, so the Ricci curvature also stays positive definite. From \eqref{kr
metric}, we obtain that
\begin{equation}\label{dtu0b}
\dt(\ul u_0+\ul b)_{i\jbar}\le0
\text{~and~similarly~}
\dt(\ul u_0)_{i\jbar}\le 0
\end{equation}
in the sense of matrices.  The second inequality follows by noting that $(\ul
u_0)_{i\jbar}=(U_0)_{i\jbar}$ and thus also corresponds to a solution to the
\eqref{krfequation} with positive holomorphic bisectional curvature. The chain
rule and the transformation formula for integrals imply that
$$I_p(t)=\int\limits_{\C^n}\vert\ul b\vert^p
\det(\ul a_{i\jbar})$$
as $\det(a_{i\jbar})=e^{-nt}\det(\ul a_{i\jbar})$
and the volume elements differ by a factor $e^{nt}$.
Here $(\ul a_{i\jbar})$ is the inverse of
$$\int\limits_0^1\left((\ul u_0)_{\cdot\cdot}+\tau
\ul b_{\cdot\cdot}\right)^{i\jbar}d\tau.$$
Note that
$$(\ul u_0)_{i\jbar}+\tau\ul b_{i\jbar}=
\tau((\ul u_0)_{i\jbar}+\ul b_{i\jbar})
+(1-\tau)(\ul u_0)_{i\jbar}$$
and we get from \eqref{dtu0b}
$$\dt\left((\ul u_0)_{i\jbar}
+\tau\ul b_{i\jbar}\right)\le0.$$
As $(\ul a_{i\jbar})$ is obtained by taking the inverse
of this matrix, integrating, and taking the inverse
once more, the definiteness for the time derivative
is inverted twice, so $\dt\ul a_{i\jbar}\le0$.
Finally, $(\ul a_{i\jbar})$ is positive definite, so it
follows that
\begin{equation}\label{dt det le zero}
\dt\det(\ul a_{i\jbar})\le0.
\end{equation}

It is not obvious, whether $I_p(t)$ is differentiable with
respect to $t$ or not. Therefore, we define for radii $R>0$
$$I_{p,R}(t):=\int\limits_{B_R}\vert\ul b\vert^p
\det(\ul a_{i\jbar}).$$
In order to compute $\dt I_{p,R}(t)$, we have to compute
the evolution equation for $\ul b$,
\begin{equation}\label{dot ul b}
\dt{\ul b}=\log\frac{\det((\ul u_0)_{i\jbar}
+\ul b_{i\jbar})}{\det((\ul u_0)_{i\jbar})}
=\ul a^{i\jbar}\ul b_{i\jbar}.
\end{equation}
Using \eqref{dt det le zero} and \eqref{dot ul b}
\begin{align*}
\dt I_{p,R}(t)=&\int\limits_{B_R}p\vert\ul b\vert^{p-2}\ul b
\left(\dt\ul b\right)\det(\ul a_{i\jbar})
+\int\limits_{B_R}\vert\ul b\vert^p\dt\det(a_{i\jbar})\\
\le&\int\limits_{B_R}p\vert\ul b\vert^{p-2}\ul b
\det(\ul a_{k\lbar})\ul a^{i\jbar}\ul b_{i\jbar}.
\end{align*}
To estimate further, we denote by $g$ the real metric
corresponding to $(\ul a_{i\jbar})$, see e.\ g.\ \cite{Shi},
and obtain in real coordinates
$$\dt I_{p,R}(t)\le\int\limits_{B_R}p\vert\ul b\vert^{p-2}
\ul b\Delta_g\ul b\sqrt{\det(g)}\,dx
\equiv\int\limits_{B_R}p\vert\ul b\vert^{p-2}
\ul b\Delta_g\ul b\,d\mu_g.$$
We apply the divergence theorem and use $\nu$ to denote
the exterior unit normal to $B_R$ with respect to the
metric $g$ which coincides with $\tfrac x{\vert x\vert}$
up to a positive factor
$$\dt I_{p,R}(t)\le-\int\limits_{B_R}p(p-1)\vert\ul b\vert^{p-2}
\langle\nabla\ul b,\nabla\ul b\rangle_g \,d\mu_g
+\int\limits_{\partial B_R}p\vert\ul b\vert^{p-2}
\ul b\langle\nabla\ul b,\nu\rangle_g d\,\mathcal{H}^{2n-1}_g.$$
Here we used suggestive invariant notation.
We apply Lemma \ref{monotonicity preserved} to see that
the boundary integral is non-positive and get
$I_{p,R}(t_1)\ge I_{p,R}(t_2)$ for $0\le t_1\le t_2$.
Finally, we let $R\to\infty$ and obtain the claimed inequality.
\end{proof}

\begin{cor}\label{Lp corollary}
Let $u_0$ be a C-soliton and
$b:\C^n\to\R$ the barrier constructed in Section
\ref{barrier construction}. Assume that $p\ge2$ is
chosen such that the $L^p$-norm
$$\Vert b\Vert_{L^p}:=\int\limits_{\C^n}
\vert b\vert^p\det((u_{0})_{i\jbar}+b_{i\jbar})$$
is finite for $t=0$.
Then the $L^p$-norm of $b$ stays uniformly bounded when
$b$ evolves by \KRF (\ref{mpma2})
$$\Vert b(t)\Vert_{L^p}\le c\cdot\Vert b(0)\Vert_{L^p},$$
where the constant depends only on the uniform equivalence
of the metrics $(u_0)_{i\jbar}$ and $(u_0)_{i\jbar}+b_{i\jbar}$
that is guaranteed during the evolution.
\end{cor}

\begin{proof}[Proof of Lemma \ref{convergencetozero} and
Theorem \ref{stab thm}]
Lemmata \ref{monotonicity preserved} and \ref{Lp corollary}
together with the arguments at the beginning of the section
complete the proof of Lemma \ref{convergencetozero} and thus
of Theorem \ref{stab thm}.
\end{proof}

\begin{appendix}

\section{Preserving Monotonicity}
\label{mon pres}

\begin{proof}[Proof of Lemma \ref{monotonicity preserved}]
It is clear that the rotational symmetry is preserved during
the evolution.

If $b$ is not a monotone decaying function of $|z|$ for
all $t>0$, we choose $0\le T<\infty$ maximal such that $b$
is monotone decaying in $|z|$ for $t\in[0,\,T]$. Note that
$b$ is clearly monotone on a relatively closed subset in
time. Our lemma follows if we can show that $b$ stays
monotone for a while after $T$. To simplify notation,
we note, that applying (the independently proven) Lemma
\ref{Lp lemma} to the time interval $[0,\,T]$, where $b$
is monotone, yields that $\lim\limits_{|z|\to\infty}
b(|z|,\,T)=0$. A similar argument works if we don't use
this fact, we just have to take into account the
possibly different $\inf b(\cdot, T)$.

First, we consider $b$ on $\C^n\setminus B_R(0)$ for $R\gg1$. The radius $R$
depends only on the fact, that certain coefficients in the ordinary differential
equation are not too far from the corresponding values in the asymptotic
expansion for the soliton. So the value of $R$ depends only on $b(\cdot,0)$ and
our initial soliton as the initial soliton and the perturbed soliton stay
uniformly equivalent during the evolution. We have $b(R,\,T)>0$ as otherwise the
strong maximum principle would imply $b(\cdot,\,T)\equiv0$, so
$b(\cdot,\,t)\equiv0$ for $t>T$, contradicting the maximality of $T$. Due to the
uniformly bounded geometry during the evolution, there exists $T^*>T$ such that
$b(R,\,t)-\tfrac12b(R,\,T)\ge\tfrac1c>0$ for $t\in\left[T,\,T^*\right]$.

Note that both $u_0$ and $\tilde b=u_0+b$ solve \eqref{mpma1}. We consider $u_0$
and $b$ as functions of $s=\log|z|^2$ and $t$ and use $\ol u_0$ and $\ol b$ to
indicate that. Equation \eqref{mpma1} implies that
\begin{align*}
\dt\ol u_0=&\log\ol u_0''+(n-1)\log\ol u_0'-ns+\ol u_0'\\
\intertext{and}
\dt\left(\ol u_0+\ol b\right)=&\log\left(\ol u_0''+\ol b''\right)
+(n-1)\log\left(\ol u_0'+\ol b'\right)
-ns+\ol u_0'+\ol b'.\\
\end{align*}
Considering the difference of these two evolution equations
gives
$$\dt\ol b=\int\limits_0^1\frac1{\ol u_0''+\tau\ol b''}d\tau
\cdot\ol b''+(n-1)\int\limits_0^1\frac1{\ol u_0'+\tau\ol b'}
d\tau\cdot \ol b'+\ol b'.$$
As $\ol u_0'=\phi$ in the notation of Section \ref{ODE sec},
we see that $\ol b$ fulfills a parabolic equation of the form
$$\dt\ol b=\alpha\ol b''+\beta\ol b',$$
where $\alpha$, $\beta$, $\alpha^{-1}\in
L^\infty((\log R^2,\infty))$ for $R\gg1$ fixed appropriately.
As $b(R,\,t)-\tfrac12b(R,\,T)\ge\tfrac1c>0$ for $T\le t
\le T^*$, we can extend $\alpha$, $\beta$, and $\ol b$ from
$[\log R^2,\infty]\times[T,\,T^*]$ to $\R\times[T,\,T^*]$ as in
the case with boundary in \cite{angenent} and apply the result
of this paper to see that for
$h\in\left(0,\tfrac12b(R,\,T)\right)$,
$\#\{r\ge R:b(r,\,t)=h\}=1$ for fixed $t\in(T,\,T^*]$. This
implies monotonicity for $r\ge R$.

It remains to prove that monotonicity is preserved
for $b>\tfrac12b(R,\,t)$.
Similarly as above, we can fix a radius $R_*>R$ and $T_*>T$
such that $b(R_*,\,t)<\tfrac12b(R,\,T)$ for $T\le t\le T_*$.
Fix $\epsilon>0$ and assume that for $t_0\in[T,\,T_*]$, there
exist $0\le r_1<r_2<R_*$ such that $b(r_2,\,t_0)\ge b(r_1,\,t_0)
+\epsilon$ and $t_0$ is chosen minimal with this property.
$b(r,\,t_0)$ tends to zero as $r\to\infty$. Choose $r_3>r_2$
minimal such that $b(r_3,\,t_0)=\tfrac12(b(r_1,\,t_0)
+b(r_2,\,t_0))$ and $r_0<r_1$ maximal such that $b(r_0,\,t_0)
=\tfrac12(b(r_1,\,t_0)+b(r_2,\,t_0))$ (if such an $r_0$
exists). Set
$\Omega:=B_{r_3}\setminus\ol{B}_{r_0}$ if $r_0$ with this
property exists, otherwise $\Omega:=B_{r_3}$. From our
assumptions, we get that
$$\osc(b,\,t,\Omega):=\sup\limits_{x\in\Omega}b(x,\,t)
-\inf\limits_{x\in\Omega}b(x,\,t)$$
is strictly smaller than $\epsilon$ for $T\le t<t_0$
and equals $\epsilon$ for $t=t_0$. Note that for $t$
close to $t_0$, $b$ is close to $\tfrac12(b(r_1,\,t_0)
+b(r_2,\,t_0))$ on $\partial\Omega$. So $\left.b\right|_{\partial
\Omega}$ does not ``contribute'' to the oscillation for
$t$ close to $t_0$ and we get a contradiction to the
strong maximum principle as a positive oscillation has to be
strictly decreasing in time (Huisken, see e.\ g.\ \cite{aw}).

As $\epsilon$ was arbitrary, we see that monotonicity is
preserved in $B_{R_*}(0)$ for $T\le t\le T_*$, so
monotonicity is preserved everywhere for $T\le t\le
\min\{T^*,\,T_*\}$ and our lemma follows.
\end{proof}

\section{Asymptotic Soliton Behavior}\label{ODE sec}

\begin{lemma}\label{phi asympt}
A solution $\phi:\R\to\R$, $\phi=\phi(s)$, fulfilling
\begin{equation}\label{phi ODE}
\phi^{n-1}\phi'e^\phi=e^{ns}
\end{equation}
and $\phi\to0$ for $s\to-\infty$ has the following
asymptotic behavior at infinity
\begin{align}
\begin{split}\label{phi C0}
\phi=&ns-{\LOG1}+\nm\frac{\LOG1}{ns}+\nm\frac1{ns}\\
&+\frac12\nm\frac{\LOG2}{n^2s^2}
-\nm(n-2)\frac{\LOG1}{n^2s^2}\\
&-\frac12\nm(3n-5)\frac1{n^2s^2}
+\frac13\nm\frac{\LOG3}{(ns)^3}\\
&-\frac12\nm(3n-5)\frac{\LOG2}{(ns)^3}\\
&+\nm\left(n^2-6n+7\right)\frac{\LOG1}{(ns)^3}\\
&+\frac16\nm\left(11n^2-46n+47\right)\frac{1}{(ns)^3}
+o\left(\frac1{s^3}\right),
\end{split}\displaybreak[1]\\
\begin{split}\label{phi C1}
\phi'=&n-\frac{n-1}s-\nm\frac{\LOG1}{ns^2}
+\nm(n-2)\frac1{ns^2}\\
&-\nm\frac{\LOG2}{n^2s^3}+\nm(3n-5)\frac{\LOG1}{n^2s^3}\\
&-\nm\left(n^2-6n+7\right)\frac1{n^2s^3}
+o\left(\frac1{s^3}\right),
\end{split}\displaybreak[1]\\
\label{phi C2}
\phi''=&\frac{n-1}{s^2}
+2\nm\frac{\LOG1}{ns^3}
-\nm(3n-5)\frac1{ns^3}
+o\left(\frac1{s^3}\right),\\
\intertext{and}
\label{phi C3}
\phi'''=&-2\frac{n-1}{s^3}+o\left(\frac1{s^3}\right).
\end{align}
\end{lemma}

We wish to emphasize that for the application we have in
mind, we don't need the high precision of \eqref{phi C0}
explicitly. But as we are not
only aiming for the asymptotic expansion for $\phi$, but
also for $\phi'$, $\phi''$ and $\phi'''$, we have to
compute the expansion for $\phi$ with high precision,
as we have to use \eqref{phi ODE} and derivatives of
this equation to determine derivatives of $\phi$
iteratively. Obviously, derivatives of the expansion
of a function do not necessarily have to coincide
with expansions of the derivatives. In our situation,
however, these two operations commute. This is
essentially due to the fact that $\phi$
satisfies \eqref{phi ODE}.

\begin{proof}
We start as in \cite{Cao}.
Separation of variables,
integration by parts and induction give
\begin{equation}\label{phi implicit}
\sum\limits_{k=0}^{n-1}(-1)^{n-k-1}\frac{n!}{k!}
\phi^ke^\phi=e^{ns}+(-1)^{n-1}n!,
\end{equation}
where the constant on the right-hand
side is chosen such that $\phi(s)\to 0$ for
$s\to-\infty$. From this formula, Cao
deduces that
$$\phi(s)=ns+o(s)\quad\text{and}\quad
\phi'(s)=n+o(1)\quad\text{for~}s\to\infty.$$

To get the asymptotic behavior of $\phi$ in \eqref{phi C0},
we can directly plug an appropriate ansatz for $\phi$ in
\eqref{phi implicit} and obtain an expression for the
next correction. This results in carrying out long
computations with increasing precision.

To verify that the expansion \eqref{phi C0} is
correct, it is convenient to rewrite
(\ref{phi implicit}) as
$$1+(-1)^{n-1}n!e^{-ns}=\left(\sum\limits_{k=0}^{n-1}
(-1)^{n-k-1}\frac{n!}{k!}\frac{\phi^k}{n(ns)^{n-1}}\right)
\left(n(ns)^{n-1}e^{\phi-ns}\right).$$
We note that \eqref{phi C0} implies
\begin{align*}
e^{\phi-ns}&n(ns)^{n-1}=1+\nm\frac{\LOG1}{ns}
+\frac{n-1}{ns}\displaybreak[0]\\
&+\frac12n\nm\frac{\LOG2}{n^2s^2}
+\nm\frac{\LOG1}{n^2s^2}\displaybreak[0]\\
&-\nm(n-2)\frac1{n^2s^2}
+\frac16(n-1)n(n+1)\frac{\LOG3}{n^3s^3}\displaybreak[0]\\
&-\frac12\nm\left(n^2-2n-1\right)\frac{\LOG2}
{n^3s^3}\displaybreak[0]\\
&-\nm\left(n^2-3\right)\frac{\LOG1}{n^3s^3}\displaybreak[0]\\
&+\frac12\nm\left(n^2-8n+11\right)\frac1{n^3s^3}
+o\left(\frac1{s^3}\right)
\end{align*}
and
\begin{align}\label{lem sum}
\begin{split}
\sum\limits_{k=0}^{n-1}&(-1)^{n-k-1}\frac{n!}{k!}
\frac{\phi^k}{n(ns)^{n-1}}=1-\nm\frac{\LOG1}{ns}-\frac{n-1}{ns}\\
&+\frac12\nm(n-2)\frac{\LOG2}{n^2s^2}
+\nm(2n-3)\frac{\LOG1}{n^2s^2}\\
&+\nm(2n-3)\frac{1}{n^2s^2}
-\frac16\nm(n-2)(n-3)\frac{\LOG3}{n^3s^3}\\
&-\frac12\nm\left(3n^2-12n+11\right)\frac{\LOG2}{n^3s^3}\\
&-2\nm(n-2)(2n-3)\frac{\LOG1}{n^3s^3}\\
&-\frac12\nm\left(7n^2-24n+21\right)\frac1{n^3s^3}
+o\left(\frac1{s^3}\right).
\end{split}
\end{align}
Moreover, it is not too complicated to see that additional
terms don't improve the approximation unless they belong
to the class $o\left(s^{-3}\right)$. Thus \eqref{phi C0}
follows.

Note that the right-hand side of (\ref{lem sum}) can also be
used for the expansion of $\frac{\exp(ns-\phi)}{n(ns)^{n-1}}$,
because
$$\sum\limits_{k=0}^{n-1}(-1)^{n-k-1}\frac{n!}{k!}
\frac{\phi^k}{n(ns)^{n-1}}=\frac{\exp(ns-\phi)}{n(ns)^{n-1}}
+o\left(\frac1{s^3}\right).$$

To determine the behavior of $\phi'$ at infinity, we note that
direct calculations give
\begin{align*}
\left(\frac{ns}\phi\right)^{n-1}=&1+\nm\frac{\LOG1}{ns}
+\frac12\nm n\frac{\LOG2}{n^2s^2}\displaybreak[0]\\
&-\nm^2\frac{\LOG1}{n^2s^2}-\nm^2\frac1{n^2s^2}\displaybreak[0]\\
&+\frac16\nm n(n+1)\frac{\LOG3}{n^3s^3}\displaybreak[0]\\
&-\frac12\nm^2(2n+1)\frac{\LOG2}{n^3s^3}\displaybreak[0]\\
&-2\nm^2\frac{\LOG1}{n^3s^3}+\frac12\nm^2(3n-5)\frac1{n^3s^3}
+o\left(\frac1{s^3}\right).
\end{align*}
Combining this with \eqref{phi ODE},
(\ref{lem sum}), and the remark following \eqref{lem sum}
gives \eqref{phi C1}.

To obtain \eqref{phi C2} and \eqref{phi C3}, we make use
of the Taylor expansion of $\frac{ns}\phi$. We get
\begin{align*}
\frac1\phi=&\frac1{ns}+\frac{\LOG1}{n^2s^2}
+\frac{\LOG2}{n^3s^3}-\nm\frac{\LOG1}{n^3s^3}\\
  &-\nm\frac1{n^3s^3}+o\left(\frac1{s^3}\right),\\
\phi''=&\phi'\left(n-\phi'-\nm\frac{\phi'}\phi\right),\\
\frac{\phi'}{\phi}=&\frac1s+\frac{\LOG1}{ns^2}
-\nm\frac1{ns^2}+\frac{\LOG2}{n^2s^3}\\
&-3\nm\frac{\LOG1}{n^2s^3}
+\nm(n-3)\frac1{n^2s^3}+o\left(\frac1{s^3}\right),\\
\phi'''&=\phi''\left(n-2\phi'-2\nm\frac{\phi'}{\phi}\right)
+\nm\left(\frac{\phi'}\phi\right)^2\phi',\\
\end{align*}
and deduce directly \eqref{phi C2} and \eqref{phi C3}.
\end{proof}

\section{Positive Holomorphic Bisectional Curvature}\label{phi inequ}

\begin{lemma}\label{lem phi b}
For the function $\phi_b$ introduced in \eqref{phi b}, we
have
\begin{align}
\phi_b>&0,\label{5-1}\\
\phi_b'>&0,\label{5-2}\\
\phi_b-\phi_b'>&0,\label{5-3}\\
\left(\phi_b'\right)^2-\phi_b\phi_b''>&0,\label{5-4}\\
\left(\phi_b''\right)^2-\phi_b'\phi_b'''>&0\label{5-5}
\end{align}
for $R\gg1$ sufficiently large.
\end{lemma}

\begin{remark}
Before we give a proof of Lemma \ref{lem phi b}, we wish to note that it implies
that $(u_0)_{i\jbar}+b_{i\jbar}$ has positive holomorphic bisectional curvature.
This follows from the calculations in \cite{Cao}. Cao gives a proof of this lemma
for a C-soliton, so it suffices to proof it in regions where we have changed
$\phi$.

Note that the proof of Lemma \ref{lem phi b} shows also
that the metric of the barrier is uniformly
equivalent to the soliton metric.
\end{remark}

\begin{proof}
\def\psiarg{\left(\frac sR\right)}
We differentiate the definition of $\phi_b$, use
Lemma \ref{phi asympt}, and get
\begin{align*}
\phi_b(s)=&ns+o(s)\mp Ks^{-1-\alpha}\alpha(2R)^\alpha
\psi\psiarg,\displaybreak[1]\\
\phi'_b(s)=&n+o(1)\pm K(1+\alpha)s^{-2-\alpha}\alpha
(2R)^\alpha\psi\psiarg\\
&\mp Ks^{-1-\alpha}\alpha(2R)^\alpha\frac1R\psi'\psiarg,
\displaybreak[1]\\
\phi''_b(s)=&\frac{n-1}{s^2}+o\left(\frac1{s^2}\right)
\mp K(1+\alpha)(2+\alpha)s^{-3-\alpha}\alpha(2R)^\alpha\psi\psiarg\\
&\pm 2K(1+\alpha)s^{-2-\alpha}\alpha(2R)^\alpha\frac1R
\psi'\psiarg\\
&\mp Ks^{-1-\alpha}\alpha(2R)^\alpha\frac1{R^2}\psi''\psiarg,
\displaybreak[1]\\
\phi'''_b(s)=&-2\frac{n-1}{s^3}+o\left(\frac1{s^3}\right)\\
&\pm K(1+\alpha)(2+\alpha)(3+\alpha)s^{-4-\alpha}
\alpha(2R)^\alpha\psi\psiarg\\
&\mp 3K(1+\alpha)(2+\alpha)s^{-3-\alpha}\alpha(2R)^\alpha
\frac1R\psi'\psiarg\\
&\pm 3K(1+\alpha)s^{-2-\alpha}\alpha(2R)^\alpha\frac1{R^2}
\psi''\psiarg\\
&\mp Ks^{-1-\alpha}\alpha(2R)^\alpha\frac1{R^3}
\psi'''\psiarg.
\end{align*}
To get \eqref{5-1}, we study $s^{-1-\alpha}R^\alpha\psi$
in detail. When we choose $R$ sufficiently large,
$|s^{-1-\alpha}R^\alpha|$ becomes arbitrarily small for
$s\ge R$. For $s\le R$, however, $\psi\left(\tfrac sR
\right)$ vanishes. Thus \eqref{5-1} follows for $s\ge R$
when $R$ is sufficiently large and is true for $s<R$
by the calculations in \cite{Cao}.

Equations \eqref{5-2}, \eqref{5-3},
and \eqref{5-4} are proved similarly.
Note, however, that the term $s^{-1-\alpha}R^{\alpha-2}\psi''$
is estimated by choosing $R$ large, as $s^{-1-\alpha}$ decays
slower as the ``leading'' term $\frac{n-1}{s^2}$ as a
function of $s$. This works as $\psi''$ is zero outside
$R\le s\le 2R$. The same arguments can also be applied
to $\psi'''$. Thus for $\phi_b$, $\phi'_b$, $\phi''_b$,
and $\phi'''_b$, the additional terms with a factor $K$
can all be absorbed in the original error terms for
$R\gg1$ fixed sufficiently large. We wish to stress, that
the sign of $\phi'''_b$, as $s\to\infty$, is important to
get \eqref{5-5}. For this reason, we had to do all the
approximations in Section \ref{ODE sec} up to such a high
precision.
\end{proof}

\begin{remark}\label{bounded rem}
The expression for the Riemannian curvature tensor for a
radially symmetric K\"ahler potential in \cite{Cao} and
the expansions of $\phi$ and $\phi_b$ at infinity imply
$\Vert Rm\Vert\le c$ for the K\"ahler metrics corresponding
to $\phi$ and $\phi_b$, respectively. Moreover, the vector
fields $\left(z^i\right)$ and $\left(z^\jbar\right)$ have
finite length with respect to these metrics.
\end{remark}

\begin{remark}\label{ip fin}
It follows from Lemma \ref{phi asympt} that
the $L^p$-norm and the uniformly equivalent
quantity considered in Corollary \ref{Lp corollary}
and Lemma \ref{Lp lemma}, respectively, are finite
for $t=0$, if $p\ge2$ is chosen so large that
\begin{equation}\label{int fin?}
\int\limits_e^\infty(\log r)^{n-1-\alpha p}
\tfrac1rdr<\infty.
\end{equation}
Choose $p$ such that $\alpha p>n+1$, with $1>\alpha>0$
as in \eqref{phi b}. Introducing a new
variable for $\log r$, we see, that the integral
in \eqref{int fin?} is finite as
$\int_1^\infty\rho^{-2}d\rho<\infty$.
\end{remark}

\end{appendix}

\bibliographystyle{amsplain}

\begin{thebibliography}{10}
\bibitem{aw} S. J. Altschuler, L.-F. Wu: Translating surfaces of
  the non-parametric mean curvature flow with prescribed contact
  angle. Calc.\ Var.\ Partial Differential Equations
  \textbf{2} (1994), 101--111.

\bibitem{angenent} S.\ Angenent: The zero set of a solution of
  a parabolic equation. J.\ Reine Angew.\ Math.\
  \textbf{390} (1988), 79--96.

\bibitem {aubin} T. Aubin: Some nonlinear problems in Riemannian
  geometry. Springer Monographs in Mathematics.
  Springer-Verlag, Berlin, 1998. xviii+395 pp.

\bibitem {cao resc} H.-D. Cao: Limits of solutions to the
  K\"ahler-Ricci flow. J. Differential Geom.\ \textbf{45} (1997),
  257--272.

\bibitem {Cao} H.-D.\ Cao: Existence of gradient
  K\"ahler-Ricci solitons. Elliptic and parabolic methods in
  geometry (Minneapolis, MN,
  1994), 1--16, A K Peters, Wellesley, MA, 1996.

\bibitem {cao flow} H.-D. Cao: Deformation of K\"ahler metrics
  to K\"ahler-Einstein metrics on compact K\"ahler manifolds.
  Invent.\ Math.\ \textbf{81} (1985), 359--372.

\bibitem {cao chow} H.-D. Cao, B. Chow: Recent developments on
  the Ricci flow. Bull.\ Amer.\ Math.\ Soc.\ (N.S.) \textbf{36}
  (1999), 59--74.

\bibitem{CY1} S.-Y. Cheng, S.-T. Yau:
  On the existence of a complete K\"ahler metric on noncompact
  complex manifolds and the regularity of Fefferman's equation.
  Comm.\ Pure Appl.\ Math.\ \textbf{33} (1980), 507--544.

\bibitem {ac} A.\ Chau: Convergence of the K\"{a}hler-Ricci
  Flow on Non-Compact Manifolds, Ph.\ D.\ Thesis, Columbia University, (2001).

\bibitem {Ecker Huisken} K. Ecker, G. Huisken: Interior estimates
  for hypersurfaces moving by mean curvature. Invent.\ Math.\
  \textbf{105} (1991), 547--569.

\bibitem {fik} M.\ Feldman, T.\ Ilmanen, D.\ Knopf: Rotationally
  symmetric shrinking and expanding gradient K\"ahler-Ricci solitons,
  preprint, {\tt http://www.math.uiowa.edu/$\sim$dknopf/}.

\bibitem{gt} D. Gilbarg, N. S. Trudinger: Elliptic partial
  differential equations of second order. Second edition.
  Grundlehren der Mathematischen Wissenschaften,
  224. Springer-Verlag, Berlin, 1983. xiii+513 pp.

\bibitem {gki} C.\ Guenther, D.\ Knopf, J.\ Isenberg:
  Stability of the Ricci flow at Ricci-flat metrics.
  Comm.\ Anal.\ Geom.\ \textbf{10} (2002), 741--777.

\bibitem {hsu} S.-Y.\ Hsu: Large time behaviour of solutions of
  the Ricci flow equation on $R\sp 2$. Pacific J. Math.\ \textbf{197}
  (2001), 25--41.

\bibitem {ham eternal} R. S. Hamilton: Eternal solutions to the
  Ricci flow. J. Differential Geom.\ \textbf{38} (1993), 1--11.

\bibitem {ham 3} R. S. Hamilton: Three-manifolds with positive Ricci
  curvature. J. Differential Geom.\ \textbf{17} (1982), 255--306.

\bibitem {os} O.\ C.\ Schn\"urer: Translating solutions to the
  second boundary value problem for curvature flows.
  Manuscripta Math.\ \textbf{108} (2002), 319--347.

\bibitem {Shi} W.-X. Shi: Ricci flow and the uniformization
  on complete noncompact K\"ahler manifolds.
  J. Differential Geom.\ \textbf{45} (1997), 94--220.

\bibitem{TY} G. Tian, S.-T. Yau:
  Existence of K\"ahler-Einstein metrics on complete
  K\"ahler manifolds and their applications
  to algebraic geometry. Mathematical aspects of string theory
  (San Diego, Calif., 1986), 574--628, Adv. Ser. Math. Phys., 1,
  World Sci. Publishing, Singapore, 1987.

\bibitem{TY2} G. Tian, S.-T. Yau:
  Complete K\"ahler manifolds with zero Ricci curvature. I.
  J. Amer.\ Math.\ Soc.\ \textbf{3} (1990), 579--609.

\bibitem {wu} L.-F.~Wu: The Ricci flow on complete $\R\sp
  2$. Comm.\ Anal.\ Geom.\ \textbf{1} (1993), 439--472.

\bibitem {yau} S.-T. Yau: On the Ricci curvature of a compact
  K\"ahler manifold and the complex Monge-Amp\`ere equation. I.
  Comm.\ Pure Appl.\ Math.\ \textbf{31} (1978), 339--411.
\end{thebibliography}

\end{document}